\documentclass[twocolumn]{autart}  
\usepackage{graphicx}          
\usepackage{amsmath}
\usepackage{amssymb}
\usepackage{bm}
\usepackage{bbm}
\usepackage{mathrsfs}
\usepackage[sort&compress]{natbib}
\usepackage{fixltx2e}
\usepackage{xcolor}

\usepackage{hyperref}
\makeatletter
\newcommand{\onetagright}{\tagsleft@false}

\makeatother

\DeclareMathOperator{\vect}{vec}
\DeclareMathOperator{\vecth}{vech}
\DeclareMathOperator{\spanv}{span}
\DeclareMathOperator{\kernl}{ker}

\newtheorem{theorem}{Theorem}[section]
\newtheorem{lemma}{Lemma}[section]

\newtheorem{proposition}{Proposition}[section]
\newtheorem{corollary}{Corollary}[section]
\newtheorem{example}{Example}
\newtheorem{assumption}{Assumption}

\begin{document}

\begin{frontmatter}
  \title{Inverse Quadratic Optimal Control for Discrete-Time Linear Systems\thanksref{footnoteinfo}} 

  \thanks[footnoteinfo]{This work is supported by China Scholarship Council.}
  \author[KTH_M]{Han Zhang}\ead{hanzhang@kth.se},
  \author[Uppsala]{Jack Umenberger}\ead{jack.umenberger@it.uu.se},
  \author[KTH_M]{Xiaoming Hu}\ead{hu@kth.se},

  \address[KTH_M]{Department of Mathematics, KTH Royal Institute of Technology, SE-100 44, Stockholm, Sweden}     
  \address[Uppsala]{Department of Information Technology, Uppsala University, Uppsala, Sweden}      

  \begin{abstract}
  In this paper, we consider the inverse optimal control problem for the discrete-time linear quadratic regulator, over finite-time horizons.
  Given observations of the optimal trajectories, and optimal control inputs, to a linear time-invariant system, the goal is to infer the parameters that define the quadratic cost function.
  The well-posedness of the inverse optimal control problem is first justified.
  In the noiseless case, when these observations are exact, we analyze the identifiability of the problem and provide sufficient conditions for uniqueness of the solution.
  In the noisy case, when the observations are corrupted by additive zero-mean noise, we formulate the problem as an optimization problem and prove the statistical consistency of the problem later.
  The performance of the proposed method is illustrated through numerical examples.

  \end{abstract}

\end{frontmatter}

\section{Introduction}

Proposed by \cite{kalman1964linear}, inverse optimal control has found a multitude of applications \citep{mombaur2010human}, \citep{finn2016guided}, \citep{berret2016don}. 
The goal of a classical optimal control problem is to find the optimal control input as well as the optimal trajectory when the cost function, system dynamics, and initial conditions are given. 
In contrast, the objective of an inverse optimal control problem is to ``reverse engineer" the cost function, given observations of optimal trajectories or control inputs, for known system dynamics. 

This paper is concerned with inverse optimal control for the discrete-time linear quadratic regulator (LQR) over finite-time horizons, i.e., finding the parameters in the quadratic objective funtion given the discrete-time linear system dynamics and (possibly noisy) observations of the optimal trajectory or control input.

Inverse optimal control for LQR, particularly in the continuous infinite time-horizon case, has been studied by a number of authors \citep{anderson2007optimal}, \citep{jameson1973inverse}, \citep{fujii1987new}. They assume the optimal feedback gain $K$ is known exactly and focus on recovering the objective function.
It was shown in \citep{boyd1994linear} that the search for matrices $Q$ and $R$ can be formulated with linear matrix inequalities (LMI) when the feedback gain $K$ is known.
\cite{priess2015solutions} consider the discrete infinite time-horizon case with noisy observations, in which the optimal feedback gain $K$ is time-invariant. 
Their approach is to identify the feedback matrix $K$ and solve for $Q$ and $R$ similar to the method proposed in \citep{boyd1994linear}.
In the finite-time horizon case, the optimal feedback gain $K_t$ is time-variant, and such an approach is not applicable.
Furthermore, the idea of ``identify the feedback gain $K_t$, then compute the corresponding $Q$" suffers from the huge number of parameters in the identification stage, i.e., the number of $K_t$'s is proportional to the length of the time-horizon. In addition, such identification does not use the knowledge that the $K_t$'s are generated by an LQR.

The contributions of this paper are three-fold.
First, we justify the well-posedness of the inverse optimal control problem for LQR in Section \ref{sec:well-poseness}.
Second, in the noiseless case (in which observations of the optimal trajectory are exact) we provide sufficient conditions for consistent estimation of the cost function, i.e., exact recovery of the matrix $Q$, c.f., Section \ref{sec:noiseless}.
Moreover, inspired by the formulation in \citep{aswani2015inverse}, we formulate the search for $Q$ as an optimization problem in the noisy case (in which observations of the optimal trajectory as well as the control input are corrupted by additive noise). We further prove that such formulation is statistically consistent, c.f., Section \ref{sec:noisy}.
The proposed method is demonstrated via a number of simulation studies, in which a better performance than the method proposed in \citep{keshavarz2011imputing} is observed, c.f. Section \ref{sec:examples}.
Conclusions are offered in Section \ref{sec:conclusion}.

\subsection{Related work}

The topic of inverse optimal control has received considerable attention in the literature.
Much of the focus, especially in recent years, has been on systems with nonlinear dynamics (of which the LQR problem is a special-case).
The authors of \citep{johnson2013inverse} consider the case of continuous finite time-horizon. 
They analyze the optimality conditions for the optimal control problem, and propose a method to minimize the violation of these conditions.
Similar ideas are used in \citep{keshavarz2011imputing} and \citep{bertsimas2015data} for the discrete finite time-horizon case (since in this case, the problem can be also interpreted as an inverse optimization problem). 
The optimization problems proposed in the above methods are numerically tractable, 
nevertheless, it has been pointed out by \citep{aswani2015inverse} that the approaches used in \citep{keshavarz2011imputing} and \citep{bertsimas2015data} are not statistically consistent and sensitive to observation noise.
\cite{aswani2015inverse} present a statistically consistent formulation, but results in a difficult optimization problem. 
\citep{molloy2016discrete}, \citep{molloy2018finite} also consider the discrete finite time-horizon case. 
They consider the Pontryagin's Maximum Principle (PMP) for the optimal control problem and pose an optimization problem whose constraints are two of the three conditions of PMP;
they then minimize the residual of the third PMP condition. In addition, they assume the optimal control input is known exactly while in our case, the optimal control input can be corrupted by noise.
The question of identifiability, i.e. uniqueness of the solution, for this approach is also addressed therein.
In a very recent work \citep{jin2018inverse}, the authors consider the discrete-time inverse optimal control problem for nonlinear systems when some segments of the trajectories and input observations are missing.
In \citep{hatz2012estimating}, the continuous finite time-horizon case is considered. 
The authors formulate the problem as a hierarchical nonlinear constrained optimization problem and two approaches that are based on PMP and Karush-Kuhn-Tucker (KKT) conditions are proposed.
The idea is to replace the inner-layer of the hierarchical optimization problem, i.e., the original optimal control problem with PMP or KKT conditions, hence making the problem tractable.
Similarly, \citep{pauwels2016linear} and \citep{rouot2017inverse} also consider the continuous finite time-horizon case, 
in which the inverse optimal control problem is studied in the framework of Hamilton-Jacobi-Bellman (HJB) equation. 
The problem is formulated as a polynomial optimization problem, and solved by a hierarchy of semidefinite relaxations.

Though our problem can be seen as a special case of the aforementioned inverse optimal control problems for nonlinear systems, we focus on the discrete finite time-horizon set-up. We also utilize the special structure of LQR to discuss the well-posedness and the identifiability of the problem. Further, inspired by \citep{aswani2015inverse}, we are able to show the statistical consistency of the estimation, which most of the papers that consider noisy observations do not cover.

\subsection{Notation}
In the remainder of the paper, $\mathbb{S}^n_+$ denotes the cone of $n$ dimensional positive semi-definite matrices and $\mathbb{S}^n$ denotes the set of all $n$ dimensional symmetric matrices. $\|\cdot\|_F$ denotes the Frobenius norm of a matrix. $\|\cdot\|$ denotes the $l_2$ norm of a vector. It holds that $\mathbb{S}^n_+\subset\mathbb{S}^n\subset\mathbb{R}^{n\times n}$, and $\mathbb{R}^{n\times n}$ is a Hilbert space whose inner-product is defined by $(G_1|G_2)=tr(G_1^TG_2)$ and $(G_1|G_1)=\|G_1\|_F^2$. 
We denote $\otimes$ as the Kronecker product and $\vect(\cdot)$, $\vecth(\cdot)$ denotes vectorization and half-vectorization respectively. It holds that $\vect(G)=\mathscr{D}\vecth(G)$, where $G\in\mathbb{S}^n$ and $\mathscr{D}$ is the duplication matrix. It holds that $\vect(G_1G_2G_3)=(G_3^T\otimes G_1)\vect(G_2)$, where $G_1\in\mathbb{R}^{k\times l}$, $G_2\in\mathbb{R}^{l\times m}$, $G_2\in\mathbb{R}^{m\times n}$. 
For a sequence of vectors $\{x_t^{(i)}\}_{t=1:N}^{i=1:M}$, where $x_t^{(i)}\in\mathbb{R}^n$, we abbreviate the notation as $x_{1:N}^{(1:M)}$ when there is no risk of confusion. And we denote $\vect(x_{1:N}^{(1:M)})$ as $[\vect(x_{1:N}^{(1)})^T,\cdots,\vect(x_{1:N}^{(M)})^T]^T$, where $\vect(x_{1:N}^{(i)})=[x_1^{(i)T},\cdots,x_N^{(i)T}]^T$. The $i$'th row of matrix $G$ is denoted as $[G]_i$.

\section{Problem Formulation and Well-Posedness}\label{sec:well-poseness}
The ``forward" optimal LQ problem reads
\begin{align}
&\min_{x_{1:N},u_{1:N-1}} J = x_N^TSx_N+{\sum}_{t=1}^{N-1} \left(u_t^TRu_t+x_t^T Qx_t\right)\label{eq:opt_ctrl_pro_fin_t}\\
&\mbox{s.t. }x_{t+1} = Ax_t+Bu_t,\quad x_1=\bar{x},\label{eq:lti}
\end{align}
where $S,Q$ are $n$-dimensional positive semidefinite matrices, $R$ is $m$-dimensional positive definite matrix, $x_t\in\mathbb{R}^n$ and $u_t\in\mathbb{R}^m$.
The inverse optimal control problem aims to find $(S,Q,R)$ given $(A,B)$, the initial value $x_1=\bar{x}$ and (possibly noisy) observations of the optimal trajectory $x_{2:N}^*$ or control input $u_{1:N-1}^*$.
For simplicity, in this paper, we consider the case of $R=I$ and $S=0$. 
In addition, it is assumed that $(A,B)$ is controllable and $B$ has full column rank. 
Moreover, we assume that $A$ is invertible.  
To see that the assumption is reasonable, consider a discrete-time system sampled from a continuous linear system $\dot{x}=\hat{A}x+\hat{B}u$, where the sample period $\Delta t$ is small. Hence for the discretized linear system, we have $A=e^{\hat{A}\Delta t}$, $B=\int_0^{\Delta t}e^{\hat{A}\tau}\hat{B}d\tau$. It is clear that $A=e^{\hat{A}\Delta t}$ is invertible. 

Before moving on considering how to solve the inverse optimal control problem, we would like to justify the well-posedness of it. 
The fundamental question for well-posedness that remains to be anwered is that: does there exist two different $Q$'s such that they can generate the same closed-loop LQR system? If there exists two different $Q$'s that can generate the same closed-loop system matrix, then the problem is obviously ill-posed.
Now we are ready to justify the well-posedness of the inverse LQR optimal control problem.
\begin{theorem}\label{thm:well_poseness}
Given the closed-loop system matrices $A_{cl}(1:N-1)$ and $N\geq n+2$, the $Q$ that is used to generate the closed-loop system matrices is unique.
\end{theorem}
\begin{pf}
We know that 
\begin{align}
K_t=-(B^TP_{t+1}B+I)^{-1}B^TP_{t+1}A,
\label{eq:K_t}
\end{align}
where $P_{2:N}\succeq 0$ is the solution to the discrete-time Riccati Equation (DRE)
\begin{equation}
\begin{aligned}
&P_t=A^TP_{t+1}A+Q-\\
&\quad A^TP_{t+1}B(B^TP_{t+1}B+I)^{-1}B^TP_{t+1}A,t=1:N-1\\
&P_N=0.
\end{aligned}
\label{eq:discrete_riccati_eq}
\end{equation}
Now assume both $Q,Q^\prime\succeq 0$ generates the closed-loop system matrices $A_{cl}(1:N-1)$. Then there are $P_{2:N}^\prime$ together with $Q^\prime$ that satisfy the DRE \eqref{eq:discrete_riccati_eq}. Denote $Q^\prime=Q+\Delta Q$, $P_t^\prime=P_t+\Delta P_t,t=2:N$, where $\Delta Q,\Delta P_{2:N}\in\mathbb{S}^n$.  

First, it is worth noticing that note that if the closed-loop systems are the same, then the control gain matrix must be the same. This is because $A+BK_t=A+BK_t^\prime\Leftrightarrow BK_t=BK_t^\prime\Rightarrow (B^TB)K_t=(B^TB)K_t^\prime$ and since $B$ has full column rank, $B^TB$ is invertible and hence $K_t=K_t^\prime$. It follows from \eqref{eq:K_t} that 
\begin{align*}
&(B^TP_{t+1}B+I)K_t=-B^TP_{t+1}A\\
\Leftrightarrow&B^TP_{t+1}(A+BK_t)=-K_t,\quad t=1:N-1.
\end{align*}
Recall $A_{cl}(t)=A+BK_t=A-B(B^TP_{t+1}B+I)^{-1}B^TP_{t+1}A$ and $P_{t+1}\in\mathbb{S}^n_+$ for all $t=1:N-1$.
Note that $A_{cl}(t)=A+BK_t$ is invertible for all $t=1:N-1$. To see that, consider the determinant of $A_{cl}(t)$:
\begin{align*}
&\det\left(A_{cl}(t)\right)=\det(A-B(B^TP_{t+1}B+I)^{-1}B^TP_{t+1}A)\\
&=\det(I-B(B^TP_{t+1}B+I)^{-1}B^TP_{t+1})\det(A)
\end{align*}
By Sylvester's determinant theorem, it follows that
\begin{align*}
&\det\left(A_{cl}(t)\right)\\
&=\det(I-B^TP_{t+1}B(B^TP_{t+1}B+I)^{-1})\det(A)\\
&=\det\left((B^TP_{t+1}B+I)-B^TP_{t+1}B\right)\\
&\times\det\left[(B^TP_{t+1}B+I)^{-1}\right]\det(A)\\
&=\det(A)/\det(B^TP_{t+1}B+I)\neq 0,\quad t=1:N-1,
\end{align*}
since $A$ is invertible.
Therefore $B^TP_{t+1}=-K_tA_{cl}^{-1}(t),t=1:N-1$. Since $Q^\prime$ generate the same $A_{cl}(1:N-1)$, it also holds that $B^TP_{t+1}^\prime=-K_tA_{cl}^{-1}(t),t=1:N-1$ and hence 
\begin{align}
B^T\Delta P_t=0,t=2:N.\label{eq:B^T_delta_P_t}
\end{align}

Moreover, recall that $P_{2:N}^\prime$ and $Q^\prime$ satisfy the DRE and $Q^\prime=Q+\Delta Q$, $P_t^\prime=P_t+\Delta P_t$. The DRE for $Q^\prime$ and $P_{2:N}^\prime$ reads
\begin{align*}
&P_{t}+\Delta P_t=A^T(P_{t+1}+\Delta P_{t+1})A+(Q+\Delta Q)\\
&-A^T(P_{t+1}+\Delta P_{t+1})B\left[B^T(P_{t+1}+\Delta P_{t+1})B+I\right]^{-1}\\
&\times B^T(P_{t+1}+\Delta P_{t+1})A=0,\quad P_N+\Delta P_N=0.
\end{align*}
By \eqref{eq:B^T_delta_P_t}, it follows from the above equation that 
\begin{equation}
\begin{aligned}
&A^T\Delta P_{t+1}A-\Delta P_t+\Delta Q=0,\quad t=2:N-1,\\
&\Delta P_N=0.
\end{aligned}
\label{eq:necessary_cond_for_uniqueness}
\end{equation}
By examining the recursion \eqref{eq:necessary_cond_for_uniqueness}, utilizing the fact $A$ is invertible and \eqref{eq:B^T_delta_P_t}, we know that
\begin{align}
&\Delta P_{N-1}=\Delta Q, \quad B^T\Delta P_{N-1}=B^T\Delta Q=0,\label{eq:delta_Q=0(1)}\\
&\Delta P_{N-2}=A^T\Delta P_{N-1}A+\Delta Q, \nonumber\\
& B^T\Delta P_{N-2}=B^TA^T\Delta QA+B^T\Delta Q=B^TA^T\Delta QA=0,\nonumber\\
& \implies B^TA^T\Delta Q=0,
\end{align}
\begin{align}
&\Delta P_{N-3}=A^T\Delta P_{N-2}A+\Delta Q\nonumber\\
&=(A^T)^2\Delta QA^2+A^T\Delta QA+\Delta Q,\nonumber\\
&B^T\Delta P_{N-3}=B^T(A^T)^2\Delta QA^2+\underbrace{B^TA^T\Delta QA+B^T\Delta Q}_{=0}\nonumber\\
&=B^T(A^T)^2\Delta Q A^2=0\implies B^T(A^T)^2\Delta Q=0\label{eq:delta_Q=0(2)}\\
&\qquad\qquad\qquad\qquad \vdots\nonumber\\
&\Delta P_2=A^T\Delta P_3A+\Delta Q\nonumber\\
&=(A^T)^{N-3}\Delta QA^{N-3}+\cdots A^T\Delta QA+\Delta Q,\nonumber
\end{align}
\begin{align}
&B^TP_2=B^T(A^T)^{N-3}\Delta QA^{N-3}\nonumber\\
&+\underbrace{B^T\left((A^T)^{N-4}\Delta QA^{N-4}+\cdots A^T\Delta QA+\Delta Q\right)}_{=0}\nonumber\\
&=B^T(A^T)^{N-3}\Delta QA^{N-3}=0\implies B^T(A^T)^{N-3}\Delta Q=0.\label{eq:delta_Q=0(3)}
\end{align}
Stacking \eqref{eq:delta_Q=0(1)}-\eqref{eq:delta_Q=0(3)} together, we get
\begin{align*}
\underbrace{\begin{bmatrix}
B^T\\B^TA^T\\\vdots\\B^T(A^T)^{N-3}
\end{bmatrix}}_{\tilde{\Gamma}}\Delta Q=0
\end{align*}
Since $(A,B)$ is controllable and $N\ge n+2$, $\tilde{\Gamma}$ has full column rank and hence $\Delta Q=0$. Thus the statement follows.
\qed
\end{pf}
\section{Inverse Optimal Control in the Noiseless Case}\label{sec:noiseless}
After justifying the well-posedness of the inverse optimal control problem, in this section, we consider inverse optimal control for the LQR problem in the noiseless case.
It is assumed that we have knowledge of $M$ sets of optimal trajectories $\lbrace x_{1:N}^{(i)*},u_{1:N-1}^{(i)*}\rbrace_{i=1}^M$, i.e., $u_t^{(i)*}=K_tx_t^{(i)*}$, where $K_t$ is the optimal feedback gain. We omit the superscript ``star" in the remainder of this section to shorten the notation.

By PMP, if $u_{1:N-1}$ and $x_{1:N}$ are the optimal control and corresponding trajectory, then there exists adjoint variables $\lambda_{2:N}$ such that
\begin{equation}
\begin{aligned}
&\lambda_t=A^T\lambda_{t+1}+Qx_t,\:t=2:N-1,\\
&\lambda_N=0,\\
&u_t=-B^T\lambda_{t+1},\:t=1:N-1.
\label{eq:PMP_original}
\end{aligned}
\end{equation}
Note that in general, PMP only provides necessary optimality conditions for optimal control problems, nevertheless, since the optimal solution to the LQ optimal control problem \eqref{eq:opt_ctrl_pro_fin_t} is unique, PMP becomes also sufficient conditions for optimality.

Note that in this case, knowing $u_{1:N-1}^{(1:M)}$ and $x_1^{(1:M)}=\bar{x}^{(1:M)}$ is equivalent to knowing $x_{1:N}^{(1:M)}$. This is because when given an optimal trajectory $x_{1:N}^{(i)}$, its corresponding optimal control $u_{1:N-1}^{(i)}$ can be determined by $u_t^{(i)}=(B^TB)^{-1}B^T(x_{t+1}^{(i)}-Ax_t^{(i)})$ since $B$ has full column rank by assumption. On the other hand, when given the initial value $\bar{x}^{(1:M)}$ and $u_{1:N-1}^{(1:M)}$, we can use \eqref{eq:lti} to compute the optimal trajectory $x_{1:N}^{(1:M)}$. Hence we do not distiguish these two cases in the remainder of this section.

Based on \eqref{eq:PMP_original}, it is straight forward to solve the inverse optimal control problem, i.e., get the matrix $Q$ by solving the following feasibility SDP problem
\begin{equation}
\begin{aligned}
& \underset{\lambda_{2:N}^{(1:M)},Q\in\mathbb{S}^n_+}{\text{minimize}}
    & & 0
    & \text{subject to}
    & &\eqref{eq:PMP_original},\:i=1:M,
\end{aligned}
\end{equation}
with a slightly abuse of notation that ``subject to \eqref{eq:PMP_original}" actually means \eqref{eq:PMP_original} with a superscript $(i)$ on every $x_t$, $\lambda_t$ and $u_t$'s.
The objective function of the above problem can be any constant, without losing generality, here we let it be 0.

Though the problem is easy in the noiseless case, however, we would like to have a closer look at the identifiability of $Q$. Namely, given a set of noiseless optimal trajectories $x_{1:N}^{(1:M)}$, is there a unique positive semidefinite matrix  that corresponds to the given optimal trajectories? 
Now we give two sufficient conditions on the given trajectories $x_{1:N}^{(1:M)}$ that can be used to determine the uniqueness of $Q$.
\begin{proposition}\label{prop:AD_full_rank}
Define matrix
\begin{equation}
\begin{aligned}
&\mathscr{A}(x)=
\begin{bmatrix}
x_2^{(1)T}\\ \vdots\\x_2^{(M)T}
\end{bmatrix}
\otimes
\begin{bmatrix}
B^T\\0\\ \vdots\\0
\end{bmatrix}+
\begin{bmatrix}
x_3^{(1)T}\\ \vdots\\x_3^{(M)T}
\end{bmatrix}
\otimes
\begin{bmatrix}
B^TA^T\\B^T\\ \vdots\\0
\end{bmatrix}\\
&+\cdots+
\begin{bmatrix}
x_{N-1}^{(1)T}\\ \vdots\\x_{N-1}^{(M)T}
\end{bmatrix}
\otimes
\begin{bmatrix}
B^T(A^T)^{N-3}\\B^T(A^T)^{N-4}\\\vdots\\B^T
\end{bmatrix}.
\end{aligned}
\label{eq:A_mat_eq}
\end{equation}
If $M(N-2)m\geq n(n+1)/2$ and $\mathscr{A}(x)\mathscr{D}$ has full column rank, then the $Q\in\mathbb{S}^n_+$ that corresponds to the given optimal trajectories $x_{1:N}^{(1:M)}$ is unique, where $\mathscr{D}$ is the duplication matrix for $\mathbb{S}^n$.
\end{proposition}
\begin{pf}
By PMP \eqref{eq:PMP_original}, it follows that 
\begin{align*}
&\begin{bmatrix}
I &-A^T \\
  &I &\ddots\\
  &  &\ddots&-A^T\\
  &  &      &I
\end{bmatrix}
\begin{bmatrix}
\lambda_2^{(i)}\\\vdots\\\lambda_N^{(i)}
\end{bmatrix}
=
\begin{bmatrix}
Qx_2^{(i)}\\\vdots\\Qx_{N-1}^{(i)}\\0
\end{bmatrix},\\
\Leftrightarrow 
&\begin{bmatrix}
\lambda_2^{(i)}\\\vdots\\\lambda_N^{(i)}
\end{bmatrix}
=
\begin{bmatrix}
I &A^T &(A^T)^2 &\cdots &(A^T)^{N-2}\\
  &I &A^T&\cdots &(A^T)^{N-3}\\
  &  &\ddots &\ddots &\vdots\\
  &  &      &I&A^T\\
  &  &      & &I
\end{bmatrix}
\begin{bmatrix}
Qx_2^{(i)}\\\vdots\\Qx_{N-1}^{(i)}\\0
\end{bmatrix},\\
\Rightarrow &-\vect(u_{1:N-1}^{(i)})
=(I\otimes B^T)\vect(\lambda_{2:N}^{(i)})
\\
&=
\begin{bmatrix}
B^TQx_2^{(i)}+B^TA^TQx_3^{(i)}+\cdots +B^T(A^T)^{N-3}Qx_{N-1}^{(i)}\\
             B^TQx_3^{(i)}    +\cdots +B^T(A^T)^{N-4}Qx_{N-1}^{(i)}\\
             \vdots\\
             B^TQx_{N-1}^{(i)}\\
             0
\end{bmatrix}
\end{align*}
Using the property of vectorization and Kronecker product, we can rewrite the above equation as
\begin{align*}
-\vect(u_{1:N-1}^{(i)})&=(x_2^{(i)T}\otimes
\begin{bmatrix}
B^T\\0\\\vdots\\0
\end{bmatrix}
+x_3^{(i)T}\otimes
\begin{bmatrix}
B^TA\\B^T\\0\\\vdots\\0
\end{bmatrix}
+\cdots\\
&x_{N-1}^{(i)T}\otimes
\begin{bmatrix}
B^T(A^T)^{N-3}\\B^T(A^T)^{N-4}\\\vdots\\B^T\\0
\end{bmatrix})\vect(Q)
\end{align*}
Stacking all $u_{1:N-2}^{(i)}$ for $i=1:M$ and by $\vect(Q)=\mathscr{D}\vecth(Q)$, we get
\begin{align}
-\vect(u_{1:N-2}^{(1:M)})=\mathscr{A}(x)\mathscr{D}\vecth(Q).
\label{eq:vech_Q_equation}
\end{align}
Since $\mathscr{A}(x)\mathscr{D}$ has dimension $M(N-1)m\times n(n+1)/2$ and $M(N-1)m\geq n(n+1)/2$, $\vecth(Q)$ must be unique provided that $\mathscr{A}(x)\mathscr{D}$ has full column rank.
\qed
\end{pf}
Sometimes, $\mathscr{A}(x)\mathscr{D}$ does not necessarily have full column rank, but we can still get a unique $Q$ by the optimal trajectories $x_{1:N}^{(1:M)}$ available. We give the following sufficient condition:
\begin{proposition}\label{prop:uniqueness_dual}
Suppose $\mathscr{A}(x)\mathscr{D}$ does not have full column rank, $\vecth(Q^\prime)$ is a solution to \eqref{eq:vech_Q_equation} and $\spanv\{\vecth(\Delta Q_k)\}=\kernl(\mathscr{A}(x)\mathscr{D})$, where $\vecth(\Delta Q_k)$ are linearly independent. In addition, suppose $\Phi^*\in\mathbb{S}^n_+$ is an optimal solution to
\begin{equation}
\begin{aligned}
& \underset{\Phi\in\mathbb{S}^n_+}{\text{minimize}}
    & & tr(Q^\prime\Phi)\\
    & \text{subject to}
    & &tr(\Delta Q_k\Phi)=0,k=1:\dim\left(\kernl\left(\mathscr{A}(x)\mathscr{D}\right)\right),
\end{aligned}
\label{eq:vech_Q_dual}
\end{equation}
where $rank(\Phi^*)=r$, $\Phi^*=G\:diag(\sigma_1,\cdots\sigma_r,0,\cdots,0)G^T$ and $GG^T=I$.
Define
\begin{align*}
\mathscr{N}_\Phi=\left\{G
\begin{bmatrix}
0&0\\0 &W
\end{bmatrix}G^T\:|\:W\in\mathbb{S}^{n-r}
\right\}.
\end{align*}
If $\mathscr{N}_\Phi\cap \spanv\{\Delta Q_k\}=\{0\}$, then the $Q\in\mathbb{S}^n_+$ that corresponds to the given optimal trajectories $x_{1:N}^{(1:M)}$ is unique.
\end{proposition}
\begin{pf}
Denote $\eta=\dim(\kernl(\mathscr{A}(x)\mathscr{D}))$. Since $\mathscr{A}(x)\mathscr{D}$ does not have full column rank, $\vecth(Q^\prime)$ is a solution to \eqref{eq:vech_Q_equation} and  $\vecth(\Delta Q_k)$ are linearly independent and spans $\kernl(\mathscr{A}(x)\mathscr{D})$, it holds that
\begin{align*}
\mathscr{A}(x)\mathscr{D}\left(\vecth(Q^\prime)+\sum_{k=1}^{\eta}\alpha_k\vecth(\Delta Q_k)\right)=0,\forall \alpha_k.
\end{align*}
What remains to show is that there exists a unique $\{\alpha_k\}_{k=1}^\eta$, such that $Q=Q^\prime+\sum_{k=1}^\eta\alpha_k\Delta Q_k\in\mathbb{S}^n_+$. 
Consider the following SDP
\begin{equation}
\begin{aligned}
& \underset{\{\alpha_k\}}{\text{minimize}}
    & & 0\\
    & \text{subject to}
    & &Q^\prime+\sum_{k=1}^\eta\alpha_k\Delta Q_k\succeq 0,
\end{aligned}
\label{eq:vech_Q_feasible}
\end{equation}
whose dual problem is \eqref{eq:vech_Q_dual}. If $\Phi^*$ is an optimal solution to \eqref{eq:vech_Q_dual} and $\mathscr{N}_\Phi\cap \spanv\{\Delta Q_k\}=\{0\}$, then the optimal solution is non-degenerate. Hence the primal problem has a unique solution.\citep{alizadeh1997complementarity}\qed
\end{pf}
\begin{rem}
If the ``real" $Q$ is strictly positive definite, then ``the matrix $\mathscr{A}(x)\mathscr{D}$ has full column rank" also becomes a necessary condition for the identifiability of $Q$. If $\mathscr{A}(x)$ does not have full rank, then there always exists some $\Delta Q\in \kernl(\mathscr{A}(x)\mathscr{D})$ and small enough $\varepsilon$ such that $\mathscr{A}(x)\mathscr{D}(\vecth(Q)+\varepsilon \vecth(\Delta Q))=-\vect(u_{1:N-2}^{(1:M)})$ and $\vecth(Q)+\varepsilon \vecth(\Delta Q)\in\mathbb{S}^n_+$ since $Q$ is an interior point in $\mathbb{S}^n_+$.\qed
\end{rem}
\begin{example}
Here is an example that illustrates Proposition \ref{prop:uniqueness_dual}. Suppose $M=1$, $N=15$, the system matrices, the initial value and the ``real" $Q$ matrix (we denote it as $\bar{Q}$) are as follows
\begin{align*}
&A=
\begin{bmatrix}
-0.1922   &-0.2490    &1.2347\\
-0.2741   &-1.0642    &-0.2296\\
1.5301    &1.6035     &-1.5062
\end{bmatrix},
B=
\begin{bmatrix}
-0.4446\\-0.1559\\0.2761
\end{bmatrix},\\
&\bar{Q}=
\begin{bmatrix}
0.0068   &-0.0116   &-0.0102\\
-0.0116  &0.0197    &0.0174\\
-0.0102  &0.0174    &0.0154
\end{bmatrix},
x_0=
\begin{bmatrix}
-25.0136 \\ -18.9592 \\ -14.8221
\end{bmatrix}.
\end{align*}
In this case, $dim(\kernl(\mathscr{A}(x)\mathscr{D}))=1$ and
\begin{align*}
&\Delta Q=
\begin{bmatrix}
0.0723   &-0.6085   &-0.1447\\
-0.6085  &-0.0422   &-0.6661\\
-0.1447  &-0.6661   &-0.3976
\end{bmatrix},\\
&\Phi^*=
\begin{bmatrix}
7.5572    &1.6696    &3.1474\\
1.6696    &4.4056    &-3.8723\\
3.1474    &-3.8723   &6.4792
\end{bmatrix},
\end{align*}
$rank(\Phi^*)=2$. If we solve the following problem 
\begin{align*}
& \underset{\beta,W\in\mathbb{R}}{\text{minimize}}
    & & 0\\
    & \text{subject to}
    & &G\begin{bmatrix}
0&0\\0 &W
\end{bmatrix}G^T=\beta \Delta Q,
\end{align*}
we will find that the only feasible solution is $\beta=W=0$. And if one solves the inverse optimal control problem, she will get an unique solution $Q^*=\bar{Q}$.\qed
\end{example}
Note that $\mathscr{A}(x)\mathscr{D}$ depends on the data. Though it has been stated in Proposition \ref{prop:AD_full_rank} that we would have a unique $Q$ that corresponds to the given optimal trajectories $x_{1:N}^{(1:M)}$ if $\mathscr{A}(x)\mathscr{D}$ has full column rank, we would like to say a bit more about the data set $x_{1:N}^{(1:M)}$, more precisely, under what conditions of the data set $x_{1:N}^{(1:M)}$ will let $\mathscr{A}(x)\mathscr{D}$ have full column rank. 
Since $\mathscr{D}$ has full column rank, $\mathscr{A}(x)\mathscr{D}$ would have full column rank if $\mathscr{A}(x)$ has full column rank. In the following we will focus on discussing what kind of data set $x_{1:N}^{(1:M)}$ would let $\mathscr{A}(x)$ have full column rank.
Before we give a sufficient condition for that, we would like to present the following lemma:
\begin{lemma}\label{lem:struct_mat_linear_indep}
If the vectors
\begin{align*}
\begin{bmatrix}
a^{(1)}_1\\a^{(1)}_2\\\vdots\\a^{(1)}_n
\end{bmatrix},
\begin{bmatrix}
a^{(2)}_1\\a^{(2)}_2\\\vdots\\a^{(2)}_n
\end{bmatrix},
\cdots,
\begin{bmatrix}
a^{(n)}_1\\a^{(1)}_2\\\vdots\\a^{(n)}_n
\end{bmatrix}
\end{align*}
are linearly independent, then matrix
\begin{align*}
\bar{\chi}=
\begin{bmatrix}
\bar{\chi}_1^{(1)} &\bar{\chi}_2^{(1)} &\cdots &\bar{\chi}_n^{(1)}\\
\bar{\chi}_1^{(2)} &\bar{\chi}_2^{(2)} &\cdots &\bar{\chi}_n^{(2)}\\
\vdots&\vdots&\vdots\\
\bar{\chi}_1^{(n)} &\bar{\chi}_2^{(1)} &\cdots &\bar{\chi}_n^{(n)}
\end{bmatrix}
\end{align*}
is nonsingular $\forall \xi_{j,l}^{(i)}$, where
\begin{align*}
\bar{\chi}_j^{(i)}=
\begin{bmatrix}
a_j^{(i)} &\xi_{j,1}^{(i)} &\cdots &\xi_{j,n-1}^{(i)}\\
          &a_j^{(i)}       &\ddots &\vdots\\
          &                &\ddots &\xi_{j,1}^{(i)}\\
          &                &       &a_j^{(i)}
\end{bmatrix}.
\end{align*}
\end{lemma}
\begin{pf}
Suppose there exists constants $\eta_{1:n}^{(1:n)}$ such that $\sum_{i=1}^n\sum_{j=1}^n\eta_j^{(i)}\left[\bar{\chi}\right]_{(i-1)n+j}=0$.
Recall the structure of $\bar{\chi}$, it must hold for the first row of every row block $\left[\bar{\chi}_1^{(i)} \bar{\chi}_2^{(i)}\cdots\bar{\chi}_n^{(i)}\right]$ in $\bar{\chi}$ that $\sum_{i=1}^n\eta_1^{(i)}\left[\bar{\chi}_1^{(i)} \bar{\chi}_2^{(i)}\cdots\bar{\chi}_n^{(i)}\right]_1=0$, which reads precisely as
\begin{align*}
&a_1^{(1)}\eta_1^{(1)}+a_1^{(2)}\eta_1^{(2)}+\cdots a_1^{(n)}\eta_1^{(n)}=0\\
&\vdots\\
&a_n^{(1)}\eta_1^{(1)}+a_n^{(2)}\eta_1^{(2)}+\cdots a_n^{(n)}\eta_1^{(n)}=0\\
\Leftrightarrow 
&\begin{bmatrix}
a_1^{(1)} &a_1^{(2)} &\cdots &a_1^{(n)}\\
a_2^{(1)} &a_2^{(2)} &\cdots &a_2^{(n)}\\
\vdots & & &\vdots\\
a_n^{(1)} &a_n^{(2)} &\cdots &a_n^{(n)}
\end{bmatrix}
\begin{bmatrix}
\eta_1^{(1)}\\\eta_1^{(2)}\\\vdots\\\eta_1^{(n)}
\end{bmatrix}=0.
\end{align*}
Since by assumption $\left[a^{(1)T}_1,a^{(1)T}_2,\cdots,a^{(1)T}_n\right]^T$, $\cdots$ $\left[a^{(n)T}_1,a^{(1)T}_2,\cdots\\a^{(n)T}_n\right]^T$ are linearly independent, the above matrix only have unique zero solution, i.e., $\eta_1^{(i)}=0,\:\forall i$. 
This implies that $\sum_{i=2}^n\sum_{j=1}^n\eta_j^{(i)}\left[\bar{\chi}\right]_{(i-1)n+j}=0$. Similar to the argument above, we can iteratively show $\eta_j^{(i)}=0,\:\forall i,j$ and hence all of the rows in $\bar{\chi}$ are linearly independent. Therefore, $\bar{\chi}$ is nonsingular.
\qed
\end{pf}
\begin{theorem}\label{thm:end_pt_indep}
Suppose $N\ge n+2$ and $M\geq n$. If there exists $n$ linearly independent $x_{N-1}^{(1:n)}$ among all $M$ sets of data, then there exists a unique $Q$ that corresponds to the given optimal trajectories $x_{1:N}^{(1:M)}$.
\end{theorem}
\begin{pf}
Recall the definition of $\mathscr{A}(x)$ in \eqref{eq:A_mat_eq}, it can be rewritten as:
\begin{align*}
\mathscr{A}(x)=\sum_{t=2}^{N-1}X_t\otimes(\mathcal{S}_t \Gamma)=\underbrace{\left[\sum_{t=2}^{N-1}X_t\otimes \mathcal{S}_t\right]}_\chi(I_n\otimes\Gamma),
\end{align*}
where
\begin{align*}
&\Gamma=
\begin{bmatrix}
B^T(A^T)^{N-3}\\B^T(A^T)^{N-4}\\\vdots\\B^T
\end{bmatrix},
\mathcal{S}_2=
\begin{bmatrix}
0_m &\cdots &0_m &I_m\\
0_m &\cdots & &0_m\\
\vdots &\vdots& &\vdots\\
0_m &\cdots & &0_m
\end{bmatrix}
\end{align*}
\begin{align*}
&\mathcal{S}_3=
\begin{bmatrix}
0_m &\cdots &I_m &0_m\\
0_m &\cdots &0_m &I_m\\
\vdots & &\vdots &\vdots\\
0_m &\cdots &0_m &0_m
\end{bmatrix},\cdots,
\mathcal{S}_{N-2}=
\begin{bmatrix}
0_m &I_m &0_m &\cdots &0_m\\
0_m &0_m &I_m &\cdots &0_m\\
\vdots &\vdots& &\ddots &I_m\\
0_m &0_m &\cdots & &0_m
\end{bmatrix},
\end{align*}
\begin{align*}
&\mathcal{S}_{N-1}=I_{m(N-2)},\: X_t=
\begin{bmatrix}
x_t^{(1)T}\\ \vdots\\x_t^{(M)T}
\end{bmatrix}.
\end{align*}

Note due to the structure of $\mathcal{S}_t$'s, $\chi$ has the following form
\begin{align*}
\chi=
\begin{bmatrix}
\chi_1^{(1)} &\cdots &\chi_n^{(1)}\\
\vdots & &\vdots\\
\chi_1^{(M)} &\cdots &\chi_n^{(M)}
\end{bmatrix},
\end{align*}
where
\begin{align*}
\chi_j^{(i)}=
\begin{bmatrix}
x_{N-1,j}^{(i)}I_m &x_{N-2,j}^{(i)}I_m &\cdots &x_{1,j}^{(i)}I_m\\
                   &x_{N-1,j}^{(i)}I_m &\ddots &\vdots\\
                   &       &\ddots &x_{N-2,j}^{(i)}I_m\\
                   &       &       &x_{N-1,j}^{(i)}I_m
\end{bmatrix}.
\end{align*}
Note that the first $n$ row blocks in $\chi$ has exactly the same structure as $\bar{\chi}$ in Lemma \ref{lem:struct_mat_linear_indep}. Since $x_{N-1}^{(1:n)}$ are linearly independent, that is, $\left[x_{N-1,1}^{(1)},\cdots,x_{N-1,n}^{(1)}\right]^T$ $\cdots$ $\left[x_{N-1,1}^{(n)},\cdots,x_{N-1,n}^{(n)}\right]^T$ are linearly independent, we can apply Lemma \ref{lem:struct_mat_linear_indep} and conclude that the matrix formed by the first $n$ row blocks of $\chi$ is nonsingular. Thus $\chi$ has full column rank.

On the other hand, since $(A,B)$ is controllable, $\Gamma$ has full column rank and $rank(\Gamma)=n$. By the property of Kronecker product, it holds that $rank(I_n\otimes\Gamma)=rank(I_n)\times rank(\Gamma)=n^2$. Therefore, due to the fact that $\chi$ has full column rank, $rank(\mathscr{A}(x))=rank\left(\chi(I_n\otimes\Gamma)\right)=rank(I_n\otimes\Gamma)=n^2$, i.e., $\mathscr{A}(x)$ has full column rank. Hence the solution $\vecth(Q)$ to the equation \eqref{eq:vech_Q_equation} is unique.
\qed
\end{pf}
\begin{rem}
Theorem \ref{thm:end_pt_indep} indicates that if among $M$ trajectories, there exists $n$ trajectories such that the second last states of each, i.e., $x_{N-1}^{(1)},\cdots,x_{N-1}^{(n)}$ are linearly independent, then $Q$ is identifiable. The theorem provides a convenient way of checking the identifiability of $Q$.\qed
\end{rem}

\section{Inverse Optimal Control in the Noisy Case}\label{sec:noisy}
Now we turn our attention to the noisy case. Inspired by \citep{aswani2015inverse}, we first pose the inverse optimal control problem in the noisy case.
Suppose the probability space $(\Omega,\mathcal{F},\mathbb{P})$ carries independent random vectors $\bar{x}\in\mathbb{R}^n$, $\{v_t\in\mathbb{R}^n\}_{t=2}^N$ and $\{w_t\in\mathbb{R}^n\}_{t=1}^{N-1}$ distributed according to some unknown distributions. The following assumptions are made in the remainder of the paper:
\begin{assumption}\label{ass:zero_mean_fin_var}
$\mathbb{E}(\|\bar{x}\|^2)<+\infty$, $\mathbb{E}(v_t)=\mathbb{E}(w_{t-1})=0$ and $\mathbb{E}(\|v_t\|^2)<+\infty$, $\mathbb{E}(\|w_{t-1}\|^2)<+\infty,t=2:N$. 
\end{assumption}
\begin{assumption}\label{ass:neighbourhood_not_measure_zero}
$\forall \eta\in\mathbb{R}^n$, $\exists r(\eta)\in\mathbb{R}\backslash\{0\}$ such that $\mathbb{P}\left(\bar{x}\in\mathscr{B}_\varepsilon(r\eta)\right)>0,\forall \varepsilon>0$, where $\mathscr{B}_\varepsilon(r\eta)$ is the open $\varepsilon$-ball centered at $r\eta$.
\end{assumption}
Equipped with the stochastic set-up above and given that the initial value $x_1$ is actually a realization of the random vector $\bar{x}$, i.e., $x_1=\bar{x}(\omega)$, the LQR problem can actually be seen as
\begin{equation}
\min \{J(u_{1:N-1}(\omega),x_{2:N}(\omega);Q;\bar{x}(\omega))|\eqref{eq:lti},\mbox{given }\omega\in\Omega\},
\label{eq:aswani_formulation}
\end{equation}
Note that the optimal control input and trajectory $\{u_t^*\},\{x_t^*\}$ are now random vectors implicitly determined by the random variable $\bar{x}$ and the parameter $Q$.
With the formulation of the ``forward problem" \eqref{eq:aswani_formulation}, we now can pose the formulation of the inverse optimal control problem.

Suppose $\{u_t^*\}$ and $\{x_t^*\}$ are corrupted by some zero mean noise, namely, $y_t=x_t^*+v_t,t=2:N$, $\mu_t=u_t^*+w_t,t=1:N-1$.
To abbreviate the notation, we denote $Y=(y_2^T,\cdots,y_N^T)^T$, $\Upsilon=(\mu_1^T,\cdots,\mu_{N-1}^T)^T$, $\xi_{x}=(\bar{x}^T,Y^T)^T$ and $\xi_u=(\bar{x}^T,\Upsilon^T)^T$. 
In addition, we assume that the ``real" $Q$ belongs to a compact set $\bar{\mathbb{S}}^n_+(\varphi)=\{Q|Q\in\mathbb{S}^n_+,\|Q\|^2_F\leq \varphi\}$. We aim to find the $Q\in\bar{\mathbb{S}}^n_+(\varphi)$ that corresponds to the optimal trajectory $\{x_t^*\}$ and control input $\{u_t^*\}$ by using the initial value $\bar{x}$ and the noisy observations $\xi_x$ or $\xi_u$.

Given $Q$ and an initial value $\bar{x}$, the solution to \eqref{eq:aswani_formulation} is unique. 
We define the risk functions $\mathscr{R}_x:\bar{\mathbb{S}}^n_+(\varphi)\mapsto\mathbb{R}$ and $\mathscr{R}_u:\bar{\mathbb{S}}^n_+(\varphi)\mapsto\mathbb{R}$ as
\begin{equation}
\mathscr{R}^x(Q)=\mathbb{E}_{\xi_x}\left[f_x(Q;\xi_x)\right],
\label{eq:risk_function}
\end{equation}
\begin{equation}
\mathscr{R}^u(Q)=\mathbb{E}_{\xi_u}\left[f_u(Q;\xi_u)\right],
\label{eq:risk_function_u}
\end{equation}
where $f_x:\bar{\mathbb{S}}^n_+(\varphi)\times \mathbb{R}^{Nn}\mapsto\mathbb{R}$ and $f_u:\bar{\mathbb{S}}^n_+(\varphi)\times \mathbb{R}^{n+m(N-1)}\mapsto\mathbb{R}$
\begin{align}
f_x(Q;\xi_x)=\sum_{t=2}^{N}\|y_t-x_t^*(Q;\bar{x})\|^2,\label{eq:risk_function1}\\
f_u(Q;\xi_x)=\sum_{t=1}^{N-1}\|\mu_t-u_t^*(Q;\bar{x})\|^2,\label{eq:risk_function1_u}
\end{align}
and $x_{2:N}^*(Q;\bar{x})$ and $u_{1:N-1}^*(Q;\bar{x})$ are the optimal solution to \eqref{eq:aswani_formulation}. In order to solve the inverse optimal control problem, we would like to minimize the risk functions, namely, 
\begin{equation}
\min_{Q\in\bar{\mathbb{S}}^n_+(\varphi)}\mathscr{R}^x(Q)
\label{eq:risk_minimization}
\end{equation}
or
\begin{equation}
\min_{Q\in\bar{\mathbb{S}}^n_+(\varphi)}\mathscr{R}^u(Q),
\label{eq:risk_minimization_u}
\end{equation}
depending on which observations are available.
Nevertheless, since the distributions of $\bar{x}$, $v_t$ and $w_t$ are unknown, the distributions of $\xi_x$ and $\xi_u$ are also unknown. We can not solve \eqref{eq:risk_minimization} and \eqref{eq:risk_minimization_u} directly. \eqref{eq:risk_function} and \eqref{eq:risk_function1_u} in principle, however, can be approximated by
\begin{align}
\mathscr{R}_M^x(Q)=\frac{1}{M}\sum_{i=1}^M f_x(Q;\xi_x^{(i)}),\label{eq:risk_function_approx}\\
\mathscr{R}_M^u(Q)=\frac{1}{M}\sum_{i=1}^M f_u(Q;\xi_u^{(i)}),\label{eq:risk_function_approx_u}
\end{align}

where $\xi_x^{(i)}$ and $\xi_u^{(i)}$ are i.i.d. random samples. We will show the statistical consistency for the approximation later.

Recall that for discrete-time LQR's in finite-time horizon, PMP \eqref{eq:PMP_original} provides sufficient and necessary conditions for optimality, hence we can express $u_{1:N-1}^*$, $x_{2:N}^*$ using \eqref{eq:PMP_original} and the approximated risk-minimizing problem reads
\begin{equation}
\begin{aligned}
\underset{Q\in\bar{\mathbb{S}}^n_+(\varphi),x_{2:N}^{(i)},\lambda_{2:N}^{(i)}}{\text{min}\qquad}
    &  \mathscr{R}_M^x(Q)=\frac{1}{M}\sum_{i=1}^M f_x(Q;\xi^{(i)})\\
     \text{s.t.\qquad}
    & x_{t+1}^{(i)}=Ax_t^{(i)}-BB^T\lambda_{t+1}^{(i)},\quad t=2:N-1,\\
    & \lambda_t^{(i)}=A^T\lambda_{t+1}^{(i)}+Qx_t^{(i)},\:t=2:N-1,\\
    & x_2^{(i)}=A\bar{x}^{(i)}-BB^T\lambda_2^{(i)},\\
    &\lambda_N^{(i)}=0,\quad i=1:M,
\end{aligned}
\label{eq:risk_minimizing_pro_approx}
\end{equation}
We omit the ``star" in the notation to avoid the confusion with the optimizer of \eqref{eq:risk_minimizing_pro_approx}. The risk-minimization problem for $\mathscr{R}_M^u(Q)$ is omitted here for the sake of brevity. The optimizer $\left(Q_M^*(\omega),x_{2:N}^{(i)*}(\omega),\lambda_{2:N}^{(i)*}(\omega)\right)$ is defined in the sense that it optimizes \eqref{eq:risk_function_approx} (or \eqref{eq:risk_function_approx_u}) for every $\omega\in\Omega$.

\begin{theorem}\label{thm:statistical_consistency}
Suppose $Q_M^*\in\bar{\mathbb{S}}^n_+(\varphi)$, $N\geq n+2$, $\{x_{2:N}^{(i)*}\}$ and $\{\lambda_{2:N}^{(i)*}\}$ solves \eqref{eq:risk_minimizing_pro_approx}, then $Q_M^*\overset{p}{\rightarrow}\bar{Q}$ as $M\rightarrow \infty$, where $\bar{Q}$ is the true value used in the ``forward" problem \eqref{eq:aswani_formulation}.
\end{theorem}
\begin{pf}
Before moving on, we would like to take a close look at \eqref{eq:PMP_original} and the system dynamics \eqref{eq:lti}. Denote $z_t=\left(x_t^T,\lambda_t^T\right)^T,t=2:N$, then the first two constraints can be written as the following implicit dynamics
\begin{equation*}
\underbrace{\begin{bmatrix}
I &BB^T\\0 &A^T
\end{bmatrix}}_{E}
z_{t+1}=
\underbrace{\begin{bmatrix}
A &0\\-Q &I
\end{bmatrix}}_F z_t,\quad t=2:N-1.
\end{equation*}
Hence we can write \eqref{eq:PMP_original} together with \eqref{eq:lti} as the following compact form
\begin{equation}
\underbrace{\begin{bmatrix}
\tilde{E} & & &\tilde{F}\\
-F &E\\
   &\ddots &\ddots\\
   &  &-F& E
\end{bmatrix}}_{\mathscr{F}(Q)}
\underbrace{\begin{bmatrix}
z_2\\\vdots\\z_N
\end{bmatrix}}_{Z}
=
\underbrace{\begin{bmatrix}
A\bar{x}^{(i)}\\0\\\vdots\\0
\end{bmatrix}}_{b(\bar{x})},
\label{eq:PMP_matrix_equation}
\end{equation}
where
\begin{align*}
\tilde{E}=
\begin{bmatrix}
I &BB^T\\0 &0
\end{bmatrix},
\tilde{F}=
\begin{bmatrix}
0 &0\\0 &I
\end{bmatrix}.
\end{align*}
We claim that $\mathscr{F}(Q)$ is invertible for all $Q\in\mathbb{S}^n_+$. Though this fact can be proven by ``brute force", i.e., by considering its determinant using Laplace expansion, perhaps the easiest way to see this is that for an arbitrary $Q\in\mathbb{S}^n_+$, \eqref{eq:PMP_matrix_equation} is a sufficient and necessary condition for the corresponding ``forward" LQR problem. Since the ``forward" LQR problem has a unique solution, it must hold that $\mathscr{F}(Q)$ is invertible for all $Q\in\mathbb{S}^n_+$. Thus, it follows that $Z=\mathscr{F}(Q)^{-1}b(\bar{x})=\mathscr{F}(Q)^{-1}\tilde{A}\bar{x}$, where $\tilde{A}=[A^T,0,\cdots,0]^T$. Hence $f_x(Q;\xi_x)$ can be rewritten as
\begin{align*}
f_x(Q;\xi_x)=\|Y-G_xZ\|^2=\|Y-G_x\mathscr{F}(Q)^{-1}\tilde{A}\bar{x}\|^2,
\end{align*}
where $G_x=I_{N-1}\otimes \left[I_n, 0_n\right]$.

It is clear that $f_x(Q;\xi_x)$ is continuous with respect to $\xi_x$, hence it is a measurable function of $\xi_x$ at each $Q$. Further, $\mathscr{F}(Q)$ is continuous and hence $\mathscr{F}(Q)^{-1}$ is continuous. Then $f_x(Q;\xi_x)$ is also continuous with respect to $Q$. 

On the other hand, since $\mathscr{F}(Q)^{-1}$ is continuous and $Q$ lives in a compact set, then $\|\mathscr{F}(Q)^{-1}\|_F$ is bounded, i.e., $\|\mathscr{F}(Q)^{-1}\|_F\leq \bar{\varphi}$ for some finite positive $\bar{\varphi}$.
It follows that $\mathbb{E}(\|Z^*\|^2)=\mathbb{E}(\|\mathscr{F}(\bar{Q})^{-1}\tilde{A}\bar{x}\|^2)\leq \|\mathscr{F}(\bar{Q})^{-1}\|_F^2\|\tilde{A}\|_F^2\mathbb{E}(\|\bar{x}\|^2)<+\infty$, where $Z^*$ corresponds to the ``true" $\bar{Q}$.

Recall that $y_t=x_t^*+v_t$ and this implies that $Y=G_xZ^*+\zeta$, where $\zeta=[v_2^T\cdots,v_N^T]^T$. By Assumption \ref{ass:zero_mean_fin_var}, $\mathbb{E}(\|v_t\|^2)<\infty$, which implies $\mathbb{E}(\|\zeta\|^2)<+\infty$. 
Therefore, $\mathbb{E}(\|Y\|^2)=\mathbb{E}(\|G_xZ^*+\zeta\|^2)\leq 2\left(\mathbb{E}(\|G_xZ^*\|^2)+\mathbb{E}(\|\zeta\|^2)\right)\leq 2\big(\|G_x\|_F^2\mathbb{E}(\|Z^*\|^2)$ $+\mathbb{E}(\|\zeta\|^2)\big)$ $<+\infty$. 
Hence it holds that 
\begin{align*}
f_x(Q,\xi_x)&=\|Y-G_x\mathscr{F}(Q)^{-1}\tilde{A}\bar{x}\|^2\\
&\leq 2\left(\|Y\|^2+\|G_x\mathscr{F}(Q)^{-1}\tilde{A}\bar{x}\|^2\right)\\
&\leq 2\left(\|Y\|^2+\|G_x\|_F^2\|\mathscr{F}(Q)^{-1}\|_F^2\|\tilde{A}\|_F^2\|\bar{x}\|^2\right)\\
&\leq 2\left(\|Y\|^2+\bar{\varphi}^2\|G_x\|_F^2\|\tilde{A}\|_F^2\|\bar{x}\|^2\right):=d(\xi_x),
\end{align*}
and it is clear that $\mathbb{E}(d(\xi_x))<+\infty$ since $\mathbb{E}(\|Y\|^2)<+\infty$ and $\mathbb{E}(\|\bar{x}\|^2)<+\infty$.
By the analysis above, we conclude that the uniform law of large numbers \citep{jennrich1969asymptotic} applies, namely, 
\begin{equation}
\sup_{Q\in\bar{\mathbb{S}}^n_+(\varphi)}\|\frac{1}{M}\sum_{i=1}^M f_x(Q,\xi_x^{(i)})-\mathbb{E}_{\xi_x}\left(f_x(Q;\xi_x)\right)\|\overset{p}{\rightarrow} 0.
\label{eq:uniform_law_of_large_numbers}
\end{equation}
Besides \eqref{eq:uniform_law_of_large_numbers}, if we are able to show $\bar{Q}$ is the unique optimizer to \eqref{eq:risk_minimization}, then $Q_M^*\overset{p}{\rightarrow}\bar{Q}$ follows directly from Theorem 5.7 in \citep{van2000asymptotic}. 

Note that by assumption, $\bar{x}$, $\{v_t\}$ are independent, hence $x_t^*(Q;\bar{x})$ are independent of the noises $\{v_t\}$. Since $y_t=x^*_t(\bar{Q},\bar{x})+v_t$, $\mathbb{E}(v_t)=0,t=2:N$, \eqref{eq:risk_function} can be simplified as $\mathscr{R}^x(Q)=L(Q)+\sum_{t=2}^N\mathbb{E}(\|v_t\|^2)$,
where $L(Q)=\mathbb{E}\left(\sum_{t=2}^N\|x_t^*(\bar{Q},\bar{x})-x_t^*(Q,\bar{x})\|^2\right)$.
It is clear that $Q=\bar{Q}$ minimizes the risk function $\mathscr{R}^x(Q)$. What remains to show is the uniqueness. By Theorem \ref{thm:well_poseness}, if $Q\neq\bar{Q}$, then $A_{cl}(1,Q)\neq A_{cl}(1,\bar{Q})$. Hence there exists $\eta\in\mathbb{R}^n$, such that $A_{cl}(1,Q)\eta\neq A_{cl}(1,\bar{Q})\eta$. 
On the other hand, by Assumption \ref{ass:neighbourhood_not_measure_zero}, $\exists r(\eta)\neq 0$, such that $\mathbb{P}(\bar{x}\in \mathscr{B}_\varepsilon(r\eta))>0,\forall \varepsilon$. Since $r\neq 0$, $A_{cl}(1,Q)(r\eta)\neq A_{cl}(1,\bar{Q})(r\eta)$.
Further, since $A_{cl}(1,Q)\eta$ is continuous with respect to $\eta,\forall Q\in\mathbb{S}^n_+$, this implies $\exists \varepsilon_1$, such that $A_{cl}(1,Q)\bar{x}\neq A_{cl}(1,\bar{Q})\bar{x},\forall \bar{x}\in\mathscr{B}_{\varepsilon_1}(r\eta)$ and $\mathbb{P}\left(\bar{x}\in\mathbb{B}_{\varepsilon_1}(r\eta)\right)>0$. Thus $L(Q)\geq \int_{\mathscr{B}_{\varepsilon_1}(r\eta)}\|\left(A_{cl}(1,Q)-A_{cl}(1,\bar{Q})\right)\bar{x}(\omega)\|^2\mathbb{P}(d\omega)>0$. Hence $\bar{Q}$ is the unique minimizer to \eqref{eq:risk_function} and the statement follows.
\qed
\end{pf}
\begin{corollary}
Suppose $Q_M^*\in\bar{\mathbb{S}}^n_+(\varphi)$, $\{x_{2:N}^{(i)*}\}$ and $\{\lambda_{2:N}^{(i)*}\}$ solves the problem of minimizing $\mathscr{R}^u(Q)$, then $Q_M^*\overset{p}{\rightarrow}\bar{Q}$ as $M\rightarrow \infty$, where $\bar{Q}$ is the true value used in the ``forward" problem \eqref{eq:aswani_formulation}.
\end{corollary}
\begin{pf}
Similar to the proof of Theorem \ref{thm:statistical_consistency}, we can rewrite $f_u(Q;\xi_u)$ as
\begin{align*}
f_u(Q;\xi_u)=\|\Upsilon-G_uZ\|^2=\|\Upsilon-G_u\mathscr{F}(Q)^{-1}\tilde{A}\bar{x}\|^2,
\end{align*}
The first part of proof that involves uniform law of large numbers can be shown analogously to the proof of Theorem \ref{thm:statistical_consistency}. What remains to show is that $\bar{Q}$ is the unique minimizer.

Analogously, since $\bar{x}$ and $\{w_t\}$ are independent by assumption, $u_t^*(Q;\bar{x})$ are independent of the noises $w_t$. Since $\mu_t=u_t^*(\bar{Q};\bar{x})+w_t$, $\mathbb{E}(w_t)=0,t=1:N-1$, the risk function \eqref{eq:risk_function1_u} can be written as $\mathscr{R}^u(Q)=L^\prime(Q)+\sum_{t=1}^{N-1}\mathbb{E}(\|w_t\|^2)$, where $L^\prime(Q)=\mathbb{E}\left(\sum_{t=1}^{N-1}\|u_t^*(\bar{Q};\bar{x})-u_t^*(Q;\bar{x})\|^2\right)$. By Theorem \ref{thm:well_poseness}, if $Q\neq \bar{Q}$, then $A_{cl}(1,Q)\neq A_{cl}(1,\bar{Q})$. Since $B$ has full column rank and recall that $A_{cl}(1,Q)=A+BK_1(Q)$, it follows that $K_1(Q)\neq K_1(\bar{Q})$. Similar to the proof of Theorem \ref{thm:statistical_consistency} and using the fact that $u_1(Q,\bar{x})=K_1(Q)\bar{x}$, we can show that $L^\prime(Q)>0$. Therefore, $\bar{Q}$ is the unique minimizer to the risk-minimizing problem and the statement follows.
\qed
\end{pf}
Now we have shown that the solutions $Q_M^*$ to the risk-minimizing problems are statistically consistent, we start to consider how to solve the problem. When actually solving the problem, $Y^{(i)}$, $\Upsilon^{(i)}$ and $\bar{x}$ are substituted with the actual measurement of the samples. From the analysis above, we know \eqref{eq:risk_minimizing_pro_approx} can be rewritten in the compact form of $\min_{Q\in\bar{\mathbb{S}}^n_+(\varphi)} \frac{1}{M}\sum_{i=1}^M\|Y^{(i)}-G_x\mathscr{F}(Q)^{-1}\tilde{A}\bar{x}\|^2$ (the risk-minimizing problem for control-input observations follows analogously).
To solve the problems, we introduce the following convex matrix function $\hat{f}_\varepsilon:\mathbb{S}^n\mapsto\mathbb{R}$, $\hat{f}_\varepsilon(Q)=\varepsilon\ln tr(e^{Q/\varepsilon})=\varepsilon \ln[\sum_{i=1}^Ne^{\sigma_i(Q)/\varepsilon}]$ \citep{nesterov2007smoothing}, where $\sigma_i(Q)$ is the $i$'th largest eigenvalue of $Q$. It holds that $\sigma_1(Q)\leq \hat{f}_\varepsilon(Q)\leq \sigma_1(Q)+\varepsilon\ln n$. Hence when $\varepsilon$ is small, the function $\hat{f}(Q)$ approximates the largest eigenvalue of $Q$ well. On the other hand, the gradient of $\hat{f}_\varepsilon(Q)$ reads $\nabla_Q\hat{f}_\varepsilon(Q)=[\sum_{i=1}^Ne^{\sigma_i(Q)/\varepsilon}]^{-1}[\sum_{i=1}^Ne^{\sigma_i(Q)/\varepsilon}\nu_i\nu_i^T],$
where $(\sigma_i(Q),\nu_i)$ are eigen-pairs of $Q$ with $\|\nu_i\|=1,\forall i$. Note that for $\varepsilon$ small enough, the gradient only numerically depends on the eigenvectors that correspond to the largest eigenvalues \citep{nesterov2007smoothing}, which makes the gradient easy to compute.
With the set-up above, we approximate the semi-positive definite constraint $Q\in\mathbb{S}^n_+$ with $\hat{f}_\varepsilon(-Q)\leq 0$ and we can solve the optimization problems with standard nonlinear optimization solvers.
\section{Numerical Examples}\label{sec:examples}
To illustrate the performance of the estimation statistically, we consider a series of discrete-time systems sampled from continuous systems $\dot{x}=\hat{A}x+\hat{B}u$ with the sampling period $\Delta t=0.1$, where
\begin{align*}
\hat{A}=
\begin{bmatrix}
0 &1\\a_1 &a_2
\end{bmatrix},
B=
\begin{bmatrix}
0\\1
\end{bmatrix};
\end{align*}
and $a_1,a_2$ are sampled from uniform distributions on $[-3,3]$. The aim for us to generate systems like this is to unsure the controllability of the systems. We take the time horizon $N=50$. The ``real " $\bar{Q}$ is generated by letting $\bar{Q}=Q_1Q_1^T$ where each elements of $Q_1$ are sampled from the uniform distribution on $[-1,1]$. We set the feasible compact set for $Q$ as $\mathbb{S}^n_+(5)$ (we discard those randomly generated $\bar{Q}$ that does not belong to $\mathbb{S}^n_+(5)$). Each element of the initial conditions $\bar{x}^{(1:M)}$ are generated by sampling from a uniform distribution supported on $[-5,5]$. We generate 200 different sets of $(\hat{A},\hat{B},\bar{Q})$ and for each fixed $(\hat{A},\hat{B},\bar{Q})$, 200 trajectories are generated, i.e., $M=200$. 15dB and 20dB of white Gaussian noises are added to $x_{2:N}^{(1:M)}$ and $u_{1:N-1}^{(1:M)}$ respectively to get $y_{2:N}^{(1:M)}$ and $\mu_{1:N-1}^{(1:M)}$. MATLAB function \texttt{fmincon} is used to solve the risk-minimizing problem. When solving the optimization problem, we use $Q=I$ as the initial iteration values for all cases.

\begin{figure}[!htpb]
	\centering
	\includegraphics[width=0.5\textwidth]{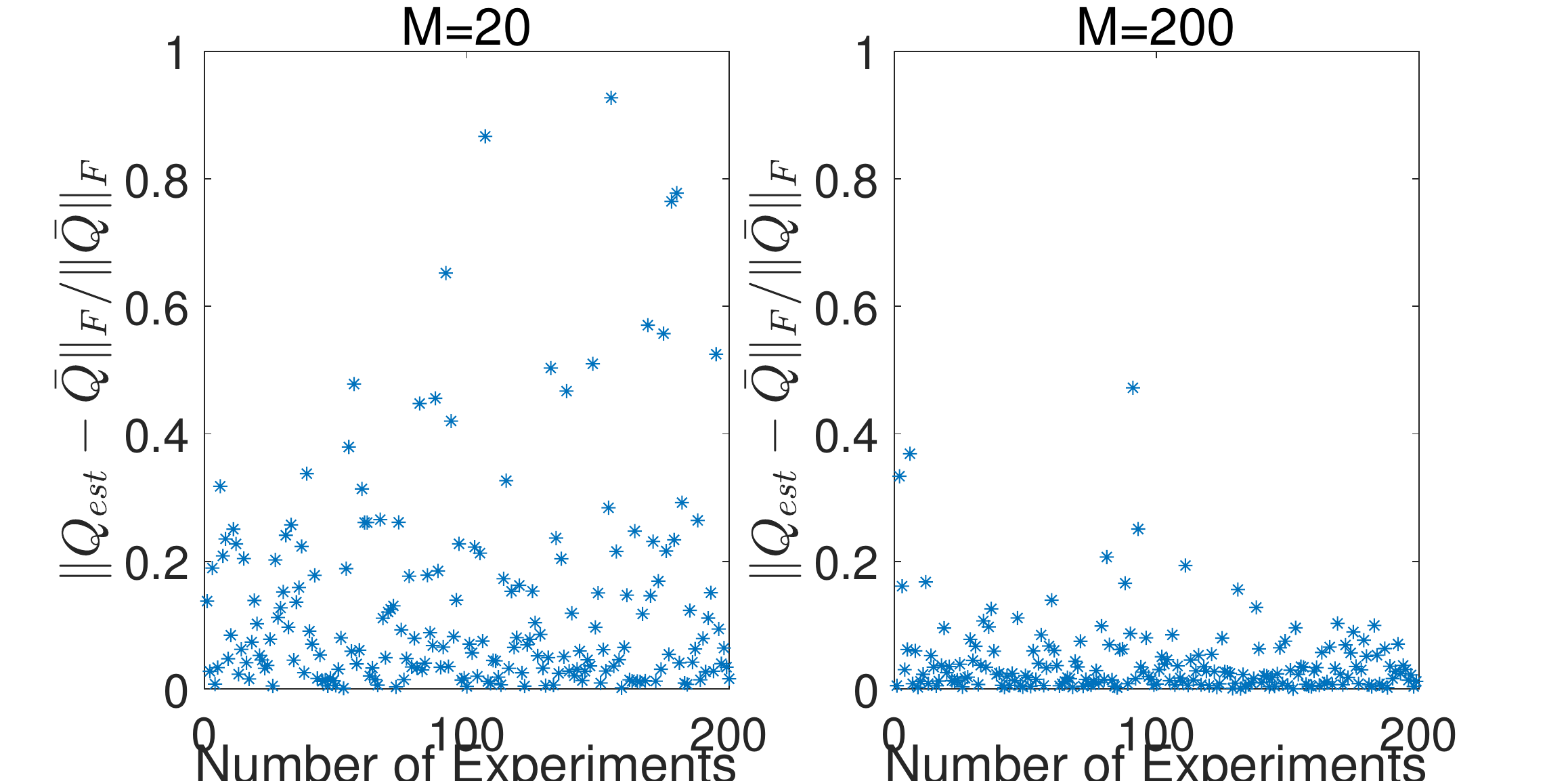}
	\caption{The relative errors of minimizing $\mathscr{R}^x_M(Q)$.}
	\label{fig:statistical_consistency1}
\end{figure}

\begin{figure}[!htpb]
	\centering
	\includegraphics[width=0.5\textwidth]{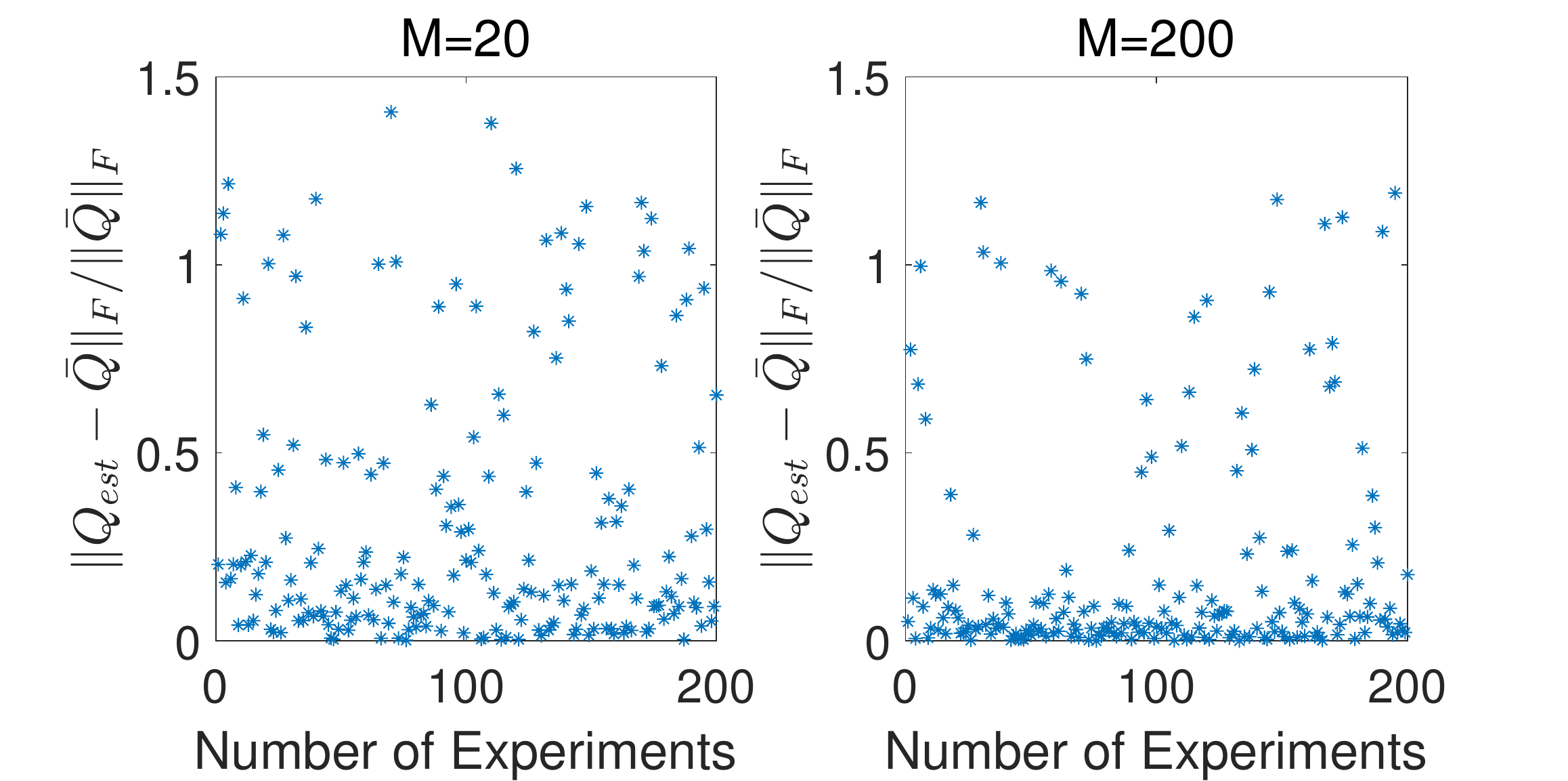}
	\caption{The relative errors of minimizing $\mathscr{R}^u_M(Q)$.}
	\label{fig:statistical_consistency2}
\end{figure}

\begin{figure}[!htpb]
	\centering
	\includegraphics[width=0.4\textwidth]{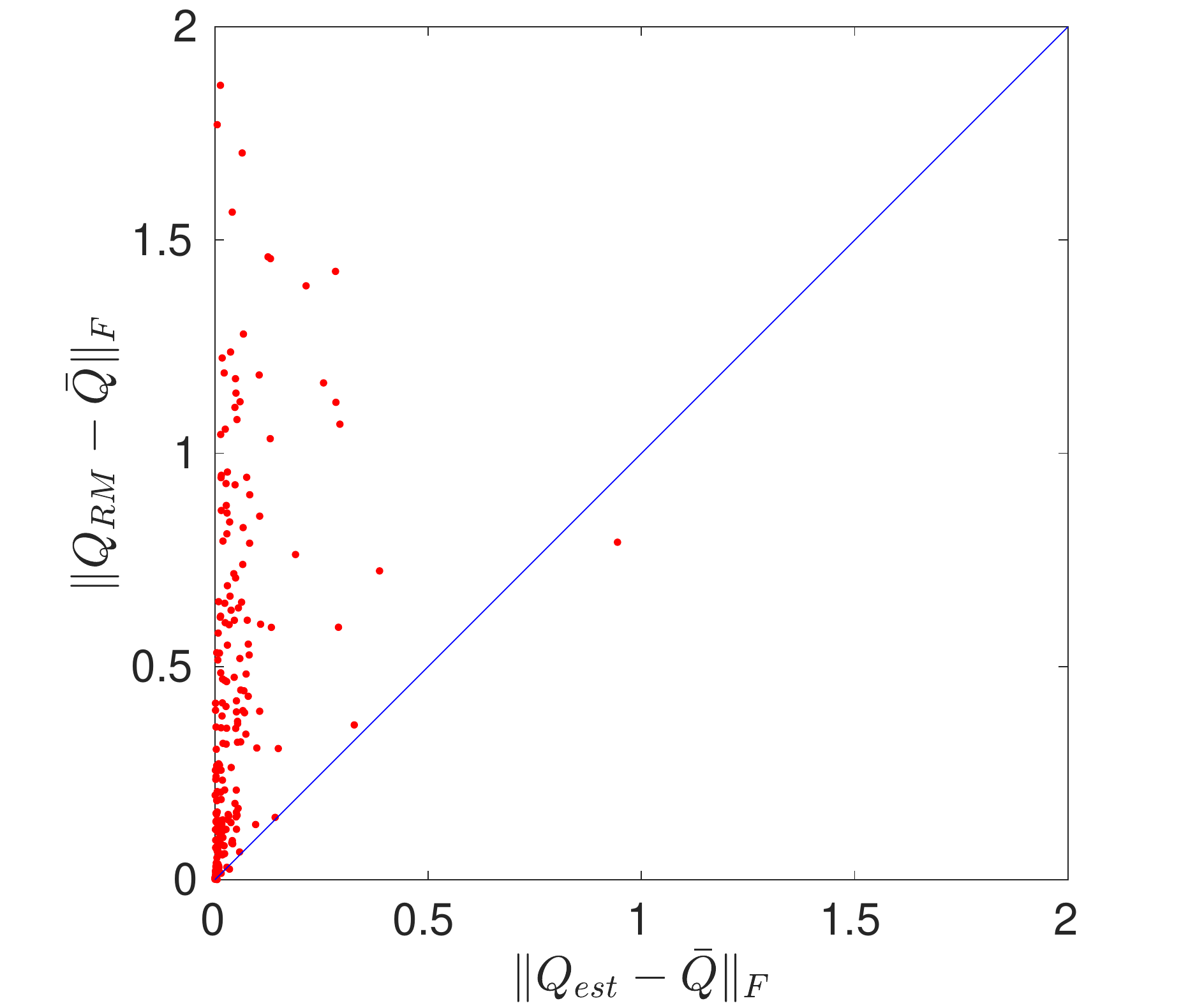}
	\caption{Our method vs. residual minimization}
	\label{fig:statistical_error1}
\end{figure}

As illustrated in Fig. \ref{fig:statistical_consistency1} and Fig. \ref{fig:statistical_consistency2}, the relative error $\|Q_{est}-\bar{Q}\|_F/\|\bar{Q}\|_F$ roughly decreases as $M$ increases. 
The result is also compared with the ``residual minimization" method proposed in \citep{keshavarz2011imputing}. In \citep{keshavarz2011imputing}, it is assumed that the observations of the solutions to the ``forward" problems are completely available, namely in this scenario, both $y_{1:N}^{(1:M)}$ and $\mu_{1:N-1}^{(1:M)}$ are available. In order to make the comparison fair, in this numerical example, observations on both of the optimal trajectories and control input are used. This will not change the statistical consistency of the method.
The result is shown in Fig. \ref{fig:statistical_error1}.
We denote the estimation of $Q$ by our method as $Q_{est}$ and the estimation by ``residual minimization" \citep{keshavarz2011imputing} as $Q_{RM}$.
In Fig. \ref{fig:statistical_error1}, the blue line illustrates $\|Q_{est}-\bar{Q}\|_F=\|Q_{RM}-\bar{Q}\|_F$. As we can see from Fig. \ref{fig:statistical_error1}, our method out-performs the residual-minimization method statistically.
\section{Conclusion}\label{sec:conclusion}
In this paper, we analyse the inverse optimal control problem for discrete-time LQR in finite-time horizons. 
We consider both the noiseless case (in which observations of the optimal trajectories are exact) and the noisy case (in which such observations are corrupted by additive noise).
The well-posedness of the problem is first justified.
In the noiseless case, we discuss identifiability of the problem, and provide sufficient conditions on the uniqueness of the solution.
In the noisy case, we formulate the search for $Q$ as an optimization problem, and prove that such formulation is statistically consistent. 
Numerical examples shows our method has a better performance than that proposed in \citep{keshavarz2011imputing}.
\bibliographystyle{model5-names}
\bibliography{ref} 
\end{document}